\documentstyle[12pt,amssymb,amsfonts]{amsart}
\newcommand{\rn}{\rr^{n}}
\newcommand{\rr}{\mbox{{\bf R}}}  

\newcommand{\norm}[1]{ \left|  #1 \right| }

\newcommand{\qtil}{{\tilde{q}}}
\newcommand{\rtil}{{\tilde{r}}}

\def\R{{\hbox{\bf R}}}

\def\rp{{r^\prime}}
\def\rpo{{r^\prime_0}}
\def\apo{{a^\prime_0}}
\def\bpo{{b^\prime_0}}
\def\rpt{{\tilde r^\prime}}

\def\qp{{q^\prime}}
\def\qpt{{\tilde q^\prime}}
\def\ap{{a^\prime}}
\def\bp{{b^\prime}}

\def\Z{{\hbox{\bf Z}}}
\def\eps{\varepsilon}
\newenvironment{proof}{\noindent {\bf Proof} }{\endprf\par}
\def \endprf{\hfill  {\vrule height6pt width6pt depth0pt}\medskip}
\def\emph#1{{\it #1}}
\def\textbf#1{{\bf #1}}

\theoremstyle{plain}
  \newtheorem{theorem}[subsection]{Theorem}
  
  \newtheorem{proposition}[subsection]{Proposition}
  \newtheorem{lemma}[subsection]{Lemma}

\theoremstyle{remark}

\theoremstyle{definition}
  \newtheorem{definition}[subsection]{Definition}

\include{psfig}

\begin{document}

\title{Low regularity semi-linear wave equations}
\author{Terence Tao}
\address{Department of Mathematics, UCLA, Los Angeles CA 90095-1555}
\email{tao@@math.ucla.edu}
\subjclass{35L70}

\begin{abstract}
We prove local well-posedness results for the semi-linear wave
equation for data in $H^\gamma$, $0 < \gamma < \frac{n-3}{2(n-1)}$,
extending the previously known results for this problem.  The improvement
comes from an introduction of a two-scale Lebesgue space $X^{r,p}_k$.
\end{abstract}
\maketitle

\section{Introduction}

We consider the initial value problem for the semi-linear wave equation
\begin{equation}
\label{semilinear}
\begin{split}
\Box u  & = F(u) \\
u(0,\cdot)  & = f \in  H^{\gamma}(\rn) \\
\partial_t u(0,\cdot) &= g \in  H^{\gamma - 1}(\rn)  
\end{split}
\end{equation}
where $n \geq 2$, $u$ is scalar or vector valued on $\R^+ \times \R^n$, 
$\Box = -\frac{\partial^2}{\partial t^2} + \Delta$ is the D'Alembertian,
$p > 1$, $\gamma \geq 0$ and the nonlinearity $F = F_p \in C^0$ satisfies\footnote{When $\gamma$ is large (e.g. $\gamma > 1/2$, or $\gamma > 3/2$) more 
regularity may be needed
on $F_p$; see \cite{lindbladsogge:semilinear}.  However, we will only be
concerned with the low-regularity problem, and such issues will not arise.}
\begin{equation}
\label{Fconditions}
F_p(0) = 0,\quad \norm{F_p(u) - F_p(v)} \lesssim \norm{u-v}(\norm{u}^{p-1}
+ \norm{v}^{p-1}).
\end{equation}

We say that the problem \eqref{semilinear} is locally well-posed in $H^\gamma$
if, for every $(f,g) \in H^\gamma \times H^{\gamma -1}$, one can find a
time\footnote{We will not concern ourselves with the exact dependence of
$T$ on the data.  In practice, one can control $T$ by the $H^\gamma \times
H^{\gamma-1}$ norm of the data unless \eqref{scaling} is satisfied with
equality, in which case $T$ depends on the data itself rather than its norm.}
$T>0$ and a unique weak solution $u \in C([0,T]; H^\gamma) \cap X$ to 
\eqref{semilinear} which depend continuously on the data, where $X$
is some additional Banach space.

The question of determining the triples $(\gamma,p,n)$ for which
\eqref{semilinear} is locally well-posed in $H^\gamma$
was studied for higher dimensions and nonlinearities by several
authors, including \cite{beals:bezard},
\cite{kapitanski:wayw}, ~\cite{lindbladsogge:semilinear}, \cite{lindblad:sharpduke}, \cite{tao:keel}. 
We summarize the known results below.

\begin{proposition}\label{known}\cite{kapitanski:wayw, lindbladsogge:semilinear, tao:keel,
lindblad:sharpduke}
 In order for \eqref{semilinear} to be locally well-posed in $H^\gamma$
for general non-linearities $F$ satisfying \eqref{Fconditions}
the following two conditions are necessary: 
\begin{align}
p(\frac{n}{2} - \gamma) &\leq \frac{n+4}{2} - \gamma
\quad\quad \text{(Scaling)}
\label{scaling}\\
p(\frac{n+1}{4} - \gamma) &\leq \frac{n+5}{4} - \gamma
\quad\quad \text{(Concentration)}
\label{concentration}
\end{align}
Conversely, if the above two conditions are satisfied and
\begin{equation}\label{ls-cond}
p(\frac{n+1}{4} - \gamma) \leq \frac{n+1}{2n} (\frac{n+3}{2} -\gamma)
\end{equation}
then (assuming sufficient regularity on $F$ if $\gamma$ is large)
\eqref{semilinear} is locally well-posed in $ H^\gamma$, with the 
exception of the case
\begin{equation}\label{hans}
 n=3, \quad p = 2, \quad \gamma = 0,
\end{equation}
which can be locally ill-posed.
\end{proposition}

For $n \geq 3$ one has the following
simultaneous endpoint of \eqref{concentration} and \eqref{ls-cond}:
\begin{equation}\label{endpoint}
\gamma  = \gamma_0 = 
\frac{n-3}{2(n-1)}\quad\quad
p  =  p_0 = \frac{(n+1)^2}{(n-1)^2 + 4}. 
\end{equation}
For $n=3$ this is \eqref{hans}, which was shown in
\cite{lindblad:sharpduke} to be locally ill-posed for
$F(u) = -|u|^2$.  For $n > 3$ 
\eqref{endpoint} was shown to be locally well-posed in \cite{tao:keel}.
The other results in the above proposition may be found in 
\cite{lindbladsogge:semilinear}, and also to a large extent in
\cite{kapitanski:wayw}.

When $n \leq 3$ or when $\gamma \geq \gamma_0$ the above results form
a complete answer to the question posed earlier, at least for
general power-type non-linearities. 
In this paper we consider the high dimension, low-regularity
case $n > 3$, $0 < \gamma < \gamma_0$.  Our main result is the following.

\begin{theorem}\label{main}  Suppose
$0 < \gamma < \gamma_0$.  Then if \eqref{concentration}
holds and
\begin{equation}\label{technical}
p(\frac{n}{4} - \gamma) \leq \frac{1}{2} (\frac{n+3}{2} -\gamma),
\end{equation}
then \eqref{semilinear} is locally well-posed in $ H^\gamma$
for all non-linearities satisfying \eqref{Fconditions},
with the possible exception of the simultaneous endpoint of 
\eqref{concentration} and \eqref{technical}
\begin{equation}\label{yuck}
\gamma 
= \frac{n+1}{4} - \frac{1}{p-1} = \frac{n+3 - \sqrt{n^2 - 2n + 33}}{8}.
\end{equation}
\end{theorem}

We note
in passing that identical results can be obtained for the semi-linear
Klein-Gordon equation by treating the mass term as an additional 
``non-linearity'', which can be treated by (e.g.) energy estimates.

\begin{figure}[htbp] \centering
\ \psfig{figure=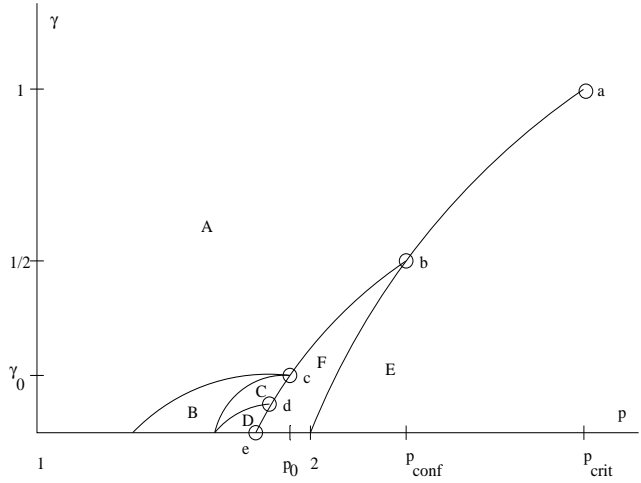,height=2.5in,width=3.3in}
\caption{Local well-posedness results for $n=4$.  }
        \label{fig:semilinearfig}
\end{figure}

These results are compared with the existing results 
in Figure \ref{fig:semilinearfig} in the case $n=4$, which is
already typical.  The scaling example (which gives \eqref{scaling})
shows that ill-posedness
is possible in the region $E$, while for non-radial 
data the concentration example (which gives \eqref{concentration}) shows
ill-posedness is possible in $F$.  (For the 
radial problem one has well-posedness everywhere above $E$; see
\cite{lindbladsogge:semilinear}).  In \cite{kapitanski:wayw} well-posedness
was shown for a certain region $A$, and extended to include $B$ in
\cite{lindbladsogge:semilinear}, including all of the boundary except
for the endpoint $c$ corresponding to \eqref{endpoint}, which was 
shown to be well-posed 
in \cite{tao:keel}.  Our results extend the positive results to the region 
$C$ including the boundary, with the exception of
the endpoint $d$ corresponding to \eqref{yuck}.  The points $a$ and $b$ 
represent the well-studied
$H^1$-critical problem and conformally invariant problem
$(\gamma,p) = (1, \frac{n+2}{n-2}), (\frac{1}{2}, \frac{n+3}{n-1})$
respectively.

We now motivate our attack strategy.  We start with the observation that
one can use standard Strichartz estimates to obtain well-posedness
for the frequency-localized equation
\begin{equation}\label{local}
 \Box u = S_j F(u), \quad u(0) = S_j f, \quad u_t(0) = S_j g
\end{equation}
all the way down to \eqref{scaling} and \eqref{concentration}; here
$S_j$ is a Littlewood-Paley projection onto a fixed frequency range
$|\xi| \sim 2^j$.  We illustrate 
this with the problem 
$$ n = 4, \quad \gamma = 0, \quad p = \frac{9}{5},$$
which is the endpoint $e$ in Figure 1.  We will use a
judiciously chosen
Strichartz estimate\footnote{The choice of exponents here is not
unique; we are using the endpoint exponents $(2,6)$ for the sake of
concreteness only.} for the linear wave equation (see \cite{tao:keel})
applied to \eqref{local}, namely
$$ \| \sqrt{-\Delta}^{-\frac{5}{54}} u \|_{L^{18}_t L^{54/25}_x} + \|u(T)\|_2
\lesssim \|\sqrt{-\Delta}^{-\frac{1}{6}} S_j F(u)\|_{L^2_t L^{6/5}_x} +
\|S_j f\|_2 + \|S_j g\|_{ H^{-1}},$$
where time is
restricted to $t \in [0,T]$ for some $T>0$.  
Because we are localizing to frequencies $|\xi| \sim 2^j$, this
estimate becomes
$$ 2^{-\frac{5}{54} j} \| u\|_{L^{18}_t L^{54/25}_x} + \|u(T)\|_2
\lesssim 2^{-\frac{1}{6} j} \| |u|^{9/5} \|_{L^2_t L^{6/5}_x} + \|f\|_2 + 
\|g\|_{ H^{-1}}.$$
Also, H\"older's inequality
gives
$$ \| |u|^{\frac{9}{5}} \|_{L^2_t L^{6/5}_x} \lesssim T^{\frac{2}{5}} 
\| u\|_{L^{18}_t L^{54/25}_x}^{\frac{9}{5}}.$$
Combining these two inequalities we obtain
$$ M + \|u(T)\|_2
\lesssim T^{\frac{2}{5}} M^{\frac{9}{5}} + \|f\|_2 + 
\|g\|_{ H^{-1}}$$
where $M=2^{-\frac{5}{54} j} \| u\|_{L^{18}_t L^{54/25}_x}$.  Thus 
a continuity argument shows that the 
$L^2$ norm of $u(T)$ is controlled by the data for sufficiently small $T$.
By adapting this inequality to differences of solutions and setting up
an iteration scheme one can also obtain local well-posedness for this
frequency-localized problem; we omit the details.

We have just seen that there are no obstructions to local well-posedness
other than
concentration and scaling if the frequencies are prevented from interacting.
To deal with the original problem \eqref{semilinear}, we must therefore
control the extent to which the $2^k$ frequency piece (say) of $F(u)$
is affected by the $2^j$ frequency piece of $u$, where $j$ is much larger
or much smaller than $k$.  Because this is a low regularity problem,
the high frequencies are less well behaved than the low frequencies,
so one expects the worst type of interaction to be when $j \gg k$.
This interaction cannot be adequately controlled by the norms used above for the
problem \eqref{local}, because of the presence of negative
derivatives.  This explains the presence of conditions such as 
\eqref{ls-cond} in previous work on the low regularity problem.

Fortunately, one can partially control this interaction with the smoothing
effect of low-frequencies.  A portion of $F(u)$ at frequency $2^k$
must necessarily be spread out at the spatial scale of $2^{-k}$, according
to the uncertainty principle.  Thus, if one takes a portion of $u$
with frequency $2^j \gg 2^k$ which is concentrated on a set which
is much ``thinner'' than $2^{-k}$, then its contribution to the
$2^k$-frequency portion of $F(u)$ will be moderated by this averaging
effect at scale $2^{-k}$.  From examining the shape of standard
examples such as the Knapp example, we see that it is indeed reasonable to 
expect the high-frequency portions of $u$ to be ``thin'', at least for the 
linear problem.

To take advantage of this effect we need
a measure of how thin the support of $u$ is compared to the spatial
scale $2^{-k}$.
To  this end we introduce a two-scale Lebesgue
space $X^{r,p}_k(\R^n)$ defined for 
$1 \leq r,p \leq \infty$ and non-negative integers $k$ by
\begin{equation}\label{x-def}
 \| u \|_{X^{r,p}_k} = \left( \sum_Q \| u \|_{L^p(Q)}^r \right)^{1/r},
\end{equation}
where $Q$ ranges over all dyadic cubes in $\R^n$ of sidelength $2^{-k}$.
(A similar norm, albeit in frequency space rather than physical space,
has appeared in \cite{borg:schrodinger},
 \cite{vargas:restrict}).  The above heuristic about the
high-frequency portion of solutions being ``thin'' can then
be captured by some Strichartz estimates for the $X^{r,p}_k$
spaces that improve upon what can be obtained by
the usual $L^r_x$ estimates and elementary inequalities.
The smoothing effect alluded to above is captured by an easy
reverse H\"older inequality
for the low-frequency pieces of functions in $X^{r,p}_k$.  These improvements
allow us to relax \eqref{ls-cond} to \eqref{technical}.

In the region $D$ in Figure $1$, \eqref{technical} fails, and the
$X^{r,p}_k$ estimates are not powerful enough to effectively
control the frequency-interference behaviour of the non-linearity.  
Indeed, it seems that one cannot go below \eqref{technical}
using norms that rely only on the size and shape of (various
frequency pieces of) $u$ and $F(u)$.  Nevertheless,
one may still conjecture that one has well-posedness in the
region $D$ (except perhaps for the endpoint $e$).  One possibility
is that the solution exhibits some additional regularity along null
directions, so that one may control it by (say) the $X^{s,b}$ spaces
as employed in \cite{beals:xsb}, \cite{borg:xsb}, \cite{kpv:xsb},
\cite{klainerman:xsb} and elsewhere; however
the non-algebraic nature of the non-linearity $F$ seems to place this
approach beyond the level of current technology, as one cannot
work exclusively in frequency space.

This paper is organized as follows.  In the next section
we set out our notation and collect many basic properties of
the $X^{r,p}_k$ spaces and the Littlewood-Paley decomposition
that we will need.  For technical reasons concerning
endpoint results we will also need a somewhat refined bilinear
interpolation theorem.  In the third section we prove
the Strichartz estimate we will need for this problem, which
involves both $X^{r,p}_k$ and $L^r_x$ spaces.  In the last section
we use this estimate together with estimates on the non-linearity to
prove the local well-posedness results.

The author wishes to thank Mark Keel, Chris Sogge, and Sergiu
Klainerman for sharing many insights about the wave equation.
This research was partly supported by NSF grant DMS-9706764 and partly
supported by MSRI (NSF grant DMS-9701955).

\section{Notation and preliminaries}

Throughout the paper, we will be working in a fixed dimension $n>3$,
and $r_0$, $\rpo$, $\gamma_0$ will denote the exponents
$$ r_0 = \frac{2(n-1)}{n-3}, \quad \rpo = \frac{2(n-1)}{n+1}, \quad
\gamma_0 = \frac{n-3}{2(n-1)}.$$
Note that $2 < r_0 < \infty$.

\begin{definition}  If $n > 3$, then an pair of exponents $(q,r)$ is called
\emph{sharp wave-admissible} if 
\begin{equation}\label{admissible}
\frac{1}{q} + \frac{(n-1)/2}{r} = 
\frac{(n-1)/2}{2}
\end{equation} and $2 \leq q, r \leq \infty$, or (equivalently)
if $(\frac{1}{q}, \frac{1}{r})$ lies on the closed line segment
between $(\frac{1}{2}, \frac{1}{r_0})$ and $(0,\frac{1}{2})$.
\end{definition}

Most of our estimates will involve sharp wave-admissible pairs
of exponents; estimates using other pairs are certainly possible,
but they can usually be obtained from the sharp estimates
via Sobolev embedding or H\"older's inequality.

For any radial function $m$, define the multiplier $m(\sqrt{-\Delta})$ by
$$ (m(\sqrt{-\Delta})f)\hat{}(\xi) = m(\xi) \hat f(\xi).$$
Define a Littlewood-Paley cutoff to be
any non-negative radial bump function supported on an annulus of
the form $\{|\xi| \sim 1\}$ which is positive on $\{\frac{1}{4} \leq |\xi| 
\leq 4\}$.
If $f$ is a function and $j$ is an integer,
we use $S_j f$ to denote the $2^j$ Littlewood-Paley frequency piece of $f$:
$$ S_j f = \beta(2^j \sqrt{-\Delta}) f;$$
for technical reasons the exact choice of $\beta$ used to define $S_j$ 
may vary from line to line, but this is not a serious problem since
a Littlewood Paley projection for one $\beta$ can always be controlled
(in virtually any space) by a finite number of such projections for any 
other $\beta$.  Henceforth we will ignore this technicality.

We also define the projection $P_0 = \phi(\sqrt{-\Delta})$, where $\phi$
is a non-negative radial bump function which equals $1$ on the ball
$\{|\xi| \leq 4\}$.

The projections $P_0$ and
$S_j$ are bounded on every $L^r_x$ space and every $X^{r,p}_k$
space, $1 \leq r, p \leq \infty$.  In particular, we have the estimate
\begin{equation}\label{triangle}
 \| f\|_{X^{r,p}_k} \lesssim \| P_0 f\|_{X^{r,p}_k} +
 \sum_{j \geq 0} \| S_j f\|_{X^{r,p}_k}
 \end{equation}
from the triangle inequality, some multiplier calculus, 
and the above observation.

We now collect some useful facts about the spaces defined in \eqref{x-def}.
Firstly, when $p=r$ these
spaces are just the Lebesgue spaces $X^{r,r}_k = L^r$.
Since $l^a \subset l^b$ for $a < b$ (by e.g. Young's
inequality) one has the inclusion 
\begin{equation}\label{inclusion}
\|f \|_{X^{b,p}_k} \lesssim \|f\|_{X^{a,p}_k} \quad \hbox{ for } a < b.
\end{equation}
By H\"older's inequality we have a similar inclusion 
for the $p$ index:
\begin{equation}\label{holder}
 \| f \|_{X^{r,p}_k} \lesssim 2^{(\frac{1}{q} - \frac{1}{p})nk} \|f\|_{X^{r,q}_k}\quad \hbox{ when } p < q.
\end{equation}
In particular, we have
\begin{equation}\label{holder-2}
 \| f \|_{X^{r,p}_k} \lesssim 2^{(\frac{1}{r} - \frac{1}{p})nk} \|f\|_r
 \quad \hbox{ when } p < r.
\end{equation}

If we localize in frequency we can reverse the above H\"older inequality
and improve\footnote{These two lemmas can also be viewed
as special cases of Sobolev embedding.} on \eqref{inclusion}.  

\begin{lemma}\label{smoothing}  (Reverse H\"older inequality)
If $1 \leq a \leq \infty$, then for any
Schwarz function $f$ and any $j \leq k$ we have
$$ \| S_j f \|_a \lesssim 2^{\frac{nk}{\ap}} \| f\|_{ X^{a,1}_k }.$$
\end{lemma}

\begin{proof}  This is trivial for $a=1$, so it suffices to verify
the case $a=\infty$.  By dilation invariance we may take $k=0$.
Since $j \leq 0$ we have the reproducing formula
$$ S_j f = S_j (f * \phi)$$
where $\phi$ is a Schwarz function whose Fourier transform equals 1
on $\{ |\xi| \leq 4\}$.  Since $S_j$ is bounded on $L^\infty$, we have
reduced ourselves to showing that
$$ \sup_x |f * \phi(x)| \lesssim \| f\|_{ X^{\infty,1}_0}.$$
Fix $x$.
  From trivial estimates we have
$$ |f * \phi(x)| \leq \sum_Q \int_Q |f(y)| |\phi(x-y)|\ dy
\lesssim \|f\|_{X^{\infty,1}_0} \sum_Q \sup_{y \in Q} |\phi(x-y)|,$$
where $Q$ ranges over unit cubes.  But from the rapid decrease of $\phi$
we have 
$$ \sum_Q \sup_{y \in Q} |\psi(x-y)| \lesssim 1$$
uniformly in $x_0$, and we are done.
\end{proof}

\begin{lemma}\label{young} (Young's inequality)
If $1 \leq a \leq b \leq \infty$, $1 \leq p \leq \infty$, and 
$k \geq 0$, then
$$
\|P_0 f \|_{X^{b,p}_k} \lesssim 2^{-nk(\frac{1}{a}-\frac{1}{b})}
\|f\|_{X^{a,p}_k}.
$$
\end{lemma}

\begin{proof}  By interpolation it suffices to prove this for $p=1$
or $p=\infty$; by duality we need only consider $p=1$.  Since
the estimate is trivial for $a=b$, we only need consider the
case $a=1$, $b=\infty$.  The
estimate now becomes
$$ \| P_0 f\|_{X^{\infty,1}_k} \lesssim 2^{-nk} \|f\|_1.$$
But this is an immediate consequence of \eqref{holder} and the
trivial estimate 
$$ \| P_0 f\|_{X^{\infty,\infty}_k} \lesssim \|f\|_1.$$
\end{proof}

Finally we observe that while the spaces $X^{r,p}_k$ are not perfectly
translation invariant, they are almost invariant in the sense that the
translation operators are uniformly bicontinuous in $X^{r,p}_k$.

We define the space-time function spaces $L^q_t L^r_x$ and $L^q_t X^{r,p}_k$
by
$$ \| F \|_{L^q_t L^r_x} = \left( \int \|F(t)\|_r^q\
dt\right)^{1/q}$$
and
$$ \| F \|_{L^q_t X^{r,2}_k} = \left( \int \|F(t)\|_{X^{r,2}_k}^q\
dt\right)^{1/q},$$
with the obvious modification for $q=\infty$.  The time integration
will usually be on a compact interval such as $0 \leq t \leq 1$.
Also we use $H^\gamma$ to denote the inhomogeneous Sobolev spaces
$(1 + \sqrt{-\Delta})^{-\gamma} L^2$, and $C(H^\gamma)$ to denote
those spacetime functions which are in $H^\gamma$ continuously
with respect to the time variable; we give $C(H^\gamma)$ the
same norm as $L^\infty_t H^\gamma$.  We will not use the
homogeneous spaces $\dot H^\gamma =  \sqrt{-\Delta}^{-\gamma} L^2$
much, although most of our results can be transferred to these spaces.

We now address the problem of interpolation between the $X^{r,p}_k$ spaces,
for fixed $k$; such interpolation was already used in the above lemmas. 
Since these spaces are equivalent to mixed Lebesgue spaces $l^r(L^p(Q))$
for a fixed $2^{-k}$-cube $Q$, the standard interpolation theorems
(e.g. the Riesz convexity theorem) apply.  In particular the spaces
$X^{r,2}_k$ behave like Hilbert-space valued $L^r$ spaces, and so obey
virtually all the interpolation identities that the scalar $L^r$ spaces do.

Finally, we will also need a certain
bilinear real interpolation theorem\footnote{It is possible to recover
the non-endpoint results in this paper without recourse to this Proposition,
or to the endpoint Strichartz estimates in \cite{tao:keel}.  More
precisely, one can prove Theorem \ref{main} using more standard 
interpolation methods provided that
\eqref{concentration} and \eqref{technical} are satisfied with
strict inequality.  We omit the details.} which we state as follows.  One
can also prove this theorem by more explicit methods; see \cite{tao:keel}.

\begin{proposition}\label{interp}  Fix $k \in \Z$ and 
$2 <  a_0, b_0 < \infty$, and 
suppose that $\{T_i(F,G): i \in \Z\}$ are a family of bilinear
forms such that one has the estimate
$$| 2^{\beta(a,b) i} T_i(F,G) | \lesssim \|F\|_{L^2_t L^{\ap}_x}
\|G\|_{L^2 X^{\bp,2}_k}$$
uniformly in $i$ for all $(\frac{1}{a}, \frac{1}{b})$ in a neighbourhood of 
$(\frac{1}{a_0}, \frac{1}{b_0})$, where $\beta(a,b)$ is an affine function of
$\frac{1}{a}$ and $\frac{1}{b}$ which is not constant with respect to 
either of the two variables.  Then one has
$$ \sum_i |2^{\beta(a_0,b_0) i} T_i(F,G)| \lesssim
 \|F\|_{L^2_t L^{\apo}_x} \|G\|_{L^2 X^{\bpo,2}_k}.$$
\end{proposition}

\begin{proof}
We introduce some notation, following \cite{bergh:interp} and
\cite{triebel:interp}.
If $A_0, A_1$ are Banach spaces contained in some larger space $A$,
we define the real interpolation
spaces $(A_0,A_1)_{\theta,q}$ for $0 < \theta < 1$,
$1 \leq q \leq \infty$ via the norm
$$ \|a\|_{(A_0,A_1)_{\theta,q}} = \bigl(\int_0^\infty (t^{-\theta} K(t,a))^q
\frac{dt}{t}\bigr)^{1/q},$$
where
$$ K(t,a) = \inf_{a = a_0 + a_1} \|a_0\|_{A_0} + t\|a_1\|_{A_1}.$$
We have the inclusions
$$ (L^2_t L^{p_0}_x, L^2_t L^{p_1}_x)_{\theta,2} = L^2_t L^{p,2}_x \subset
L^2_t L^p_x$$
whenever $p_0 \neq p_1$, $p_0, p_1 \leq 2$,
and $\frac{1}{p} = \frac{1-\theta}{p_0}
 + \frac{\theta}{p_1}$; see ~\cite{triebel:interp} Sections 1.18.2 
 and 1.18.6 for the interpolation identity, and \cite{sadosky:interp}
 for the Lorentz space inclusion.  One also has the vector-valued
 analogue of the above inclusion:
$$ (L^2_t X^{p_0,2}_k, L^2_t X^{p_1,2}_k)_{\theta,2} \subset L^2_t X^{p,2}_k.$$
 
Similarly, we have
 $$ (l^{s_0}_\infty, l^{s_1}_\infty)_{\theta,1} = l^s_1$$
 whenever $s_0 \neq s_1$ and $s = (1-\theta) s_0 + \theta
 s_1$, where
 $l^s_q = L^q(\Z, 2^{js}\ dj)$ are
 weighted sequence spaces and $dj$ is counting measure.
 See ~\cite{bergh:interp} Section 5.6.

We will use the following bilinear interpolation theorem:

\begin{lemma} (\cite{bergh:interp}, Section 3.13.5(b))
If $A_0$,$A_1$,$B_0$,$B_1$,$C_0$,$C_1$ are Banach spaces, and
the bilinear operator $T$ is bounded from
$$ T: A_0 \times B_0 \to C_0$$
$$ T: A_0 \times B_1 \to C_1$$
$$ T: A_1 \times B_0 \to C_1,$$
then one has
$$ T: (A_0,A_1)_{\theta_0,2} \times (B_0,B_1)_{\theta_1,2}
\to (C_0,C_1)_{\theta,1}$$
whenever $0 < \theta_0, \theta_1 < \theta < 1$ are
such that $\theta = \theta_0
+ \theta_1$.
\end{lemma}

Let $T(F,G)$ denote the sequence-valued bi-linear operator
$$ T(F,G) = \{ T_i(F,G) \}_{i \in \Z}.$$
Then we have
$$ T: L^2_t L^{\ap}_x \times L^2 X^{\bp,2}_k \to l^{\beta(a,b)}_\infty$$
for all $(\frac{1}{a}, \frac{1}{b})$ in a neighbourhood of
$(\frac{1}{a_0}, \frac{1}{b_0})$.
Applying the above lemma for suitable values of $(a,b)$ and using
the above inclusions, one obtains
$$ T: L^2_t L^{\ap}_x \times L^2 X^{\bp,2}_k \to l^{\beta(a,b)}_1$$
for all  $(\frac{1}{a}, \frac{1}{b})$ in a neighbourhood of
$(\frac{1}{a_0}, \frac{1}{b_0})$.  Applying this to $(a,b) = (a_0, b_0)$
one obtains the desired result.
\end{proof}

\section{Two-scale Strichartz estimates}

In this section time will always be localized to the interval
$0 \leq t \leq 1$, $f$, $g$, and $F$ will denote Schwarz functions
on $\R^n$, $\R^n$, and $[0, 1] \times \R^n$ respectively,
and $j$ and $k$ will denote \emph{non-negative} integers.

If $u$ is the solution to the linear Cauchy problem
\begin{equation}\label{standard-wave}
\Box u = F,\quad u(0) = f, \quad u_t(0) = g
\end{equation}
then we can write $u$ explicity as
\begin{equation}\label{wave-soln}
u = u_0 + \Box^{-1} F,
\end{equation}
where
\begin{align*}
u_0(t) &= \cos(t\sqrt{-\Delta}) f + \frac{\sin(t\sqrt{-\Delta})}{\sqrt{-\Delta}} g\\
\Box^{-1} F(t) &= \int_{s < t} \frac{\sin((t-s)\sqrt{-\Delta})}{\sqrt{-\Delta}} F(s)\ ds.
\end{align*}
One can localize these explicit formulae in frequency to obtain
\begin{equation}\label{pm}
\begin{split}
S_j u_0(t) &= \sum_{\pm} U^\pm_j(t) S_j f \pm i 2^{-j} U^\pm_j(t) S_j g\\
 S_j \Box^{-1} F(t) &= \sum_{\pm} \pm i 2^{-j}
\int_{s<t} U^\pm_j(t) U^\pm_j(s)^* S_j F(s)\ ds\\
\end{split}
\end{equation}
for each integer $j$, where
$$U^\pm_j(t) = \beta(2^{-j}\sqrt{-\Delta})
e^{\pm i t \sqrt{-\Delta}}$$
is a frequency localized evolution operator, and $\beta$ is a Littlewood-Paley
cutoff that varies from line to line.
Henceforth we will suppress the $\pm$ 
symbols on $U^\pm_j$.

In \cite{tao:keel} the following estimates\footnote{Strictly speaking, these
estimates was only proven (without time being localized)
in \cite{tao:keel} for $j=0$, but the general result can be recovered by 
scaling.} were proven:

\begin{proposition}\label{keel}\cite{tao:keel}
If $(q,r)$, $(\qtil, \rtil)$ are sharp wave-admissible 
pairs, and $u(t)$ is the solution to \eqref{standard-wave},
then we have the one-sided estimates
$$ 2^{-\frac{(n+1)j}{(n-1)q}} \| U_j(t) f\|_{L^q_t L^r_x} \lesssim 
\|f\|_2$$
$$ \| \int_{s<t} U_j(t) U_j(s)^* F(s)\ ds\|_{C(L^2_x)} \lesssim 
2^{\frac{(n+1)j}{(n-1)\qtil}} \|F\|_{L^\qpt_t L^\rpt_x}$$
together with the two-sided estimate
$$ 2^{-\frac{(n+1)j}{(n-1)q}} \| S_j u\|_{L^q_t L^r_x} + 
\| S_j u\|_{C(L^2_x)} \lesssim
\|S_j f\|_2 + 2^{-j} \|S_j g\|_2 + 
2^{-j} 2^{\frac{(n+1)j}{(n-1)\qtil}} \|S_j F\|_{L^\qpt_t L^\rpt_x}.$$
\end{proposition}

For examples and applications
of these estimates, see \cite{strichartz:restrictionquadratic},
 ~\cite{lindbladsogge:semilinear},~\cite{kapitanski:variablestrichartz},
 ~\cite{mss:localsmoothing},
 ~\cite{ginebre:summarywave}, ~\cite{sogge:nlwbook}, ~\cite{damiano}.

The aim of this section is to prove the analogue of this propositon
for the $X^{r,p}_k$ spaces.
We begin with the basic energy and decay estimates we will need.

\begin{lemma} If $j \geq k$, we have the energy estimate
\begin{equation}\label{energy}
 \| U_j(t) f\|_{X^{2,2}_k} \lesssim \| f\|_2,
\end{equation}
the decay estimate
\begin{equation}\label{decay}
 \| U_j(t) U_j(s)^* F\|_{X^{\infty,2}_k} \lesssim 
(2^{2k-j} |t-s|)^{-\frac{n-1}{2}} 
\|F\|_{X^{1,2}_k},
\end{equation}
and the asymmetric decay estimate
\begin{equation}\label{asymmetric}
 \| U_j(t) U_j(s)^* F\|_{X^{\infty,2}_k} \lesssim 
(2^{k-\frac{n}{n-1}j} |t-s|)^{-\frac{n-1}{2}} 
\|F\|_{1},
\end{equation}
for all $s \neq t$.
\end{lemma}

\begin{proof}
The energy estimate follows immediately from Plancherel's theorem sincee
$X^{2,2}_k = L^2$.  To prove \eqref{decay}, it suffices by self-adjointness
and interpolation
to show that
$$ \| U_j(t) U_j(s)^* F\|_{X^{\infty,1}_k} \lesssim 
(2^{2k-j} |t-s|)^{-\frac{n-1}{2}} 
\|F\|_{X^{1,1}_k}.$$
Since $X^{1,1}_k = L^1$ it suffices to verify this when $F$ is a 
delta function, which we may place at the origin since the space $X^{\infty,1}_k$ is almost translation invariant.  It now suffices to show that
$$ \| U_j(t) U_j(s)^* \delta_0 \|_{L^1(Q)} \lesssim
(2^{2k-j} |t-s|)^{-\frac{n-1}{2}}.$$ 
for all $2^{-k}$-cubes $Q$.  Similarly, \eqref{asymmetric} will follow from
$$ \| U_j(t) U_j(s)^* \delta_0 \|_{L^2(Q)} \lesssim
(2^{k-\frac{n}{n-1}j} |t-s|)^{-\frac{n-1}{2}}.$$ 
But these estimates are consequence of the standard
stationary phase estimate
$$ U_j(t) U_j(s)^* \delta_0(x) \leq C_N 2^{nj} (1+2^j |t-s|)^{-\frac{n-1}{2}}
\left(1 + 2^j( |t-s| - |x| )\right)^{-N},$$
valid for any $N > 0$.  Indeed, from these estimates we see that
$ U_j(t) U_j(s)^* \delta$ when restricted to $Q$ has a sup norm of
$O(2^{nj} (2^j |t-s|)^{-\frac{n-1}{2}})$ and is rapidly decreasing outside
a set of measure
$O(2^{-j} 2^{-(n-1)k})$, and the claimed estimates follow
from some algebraic manipulation.
\end{proof}

The estimates \eqref{energy} and \eqref{decay} imply the following family of
one-sided Strichartz estimates:
\begin{proposition}\label{one} If $(q,r)$ is a sharp wave-admissible pair and
$j \geq k$, then
\begin{equation}\label{strich-1}
 \| U_j(t) f\|_{L^q_t X^{r,2}_k} \lesssim 2^{-\frac{2k-j}{q}} \|f\|_2
\end{equation}
for all test functions $f$ on $\R^n$.
\end{proposition}

\begin{proof}  This is a special case of the abstract interpolation
theorem \cite{tao:keel}, Theorem 10.1, although strictly speaking one
must first rescale the time variable by $2^{2k-j}$ to satisfy the
conditions of that theorem.  Here we will present the
proof for $q > 2$, so that we are excluding the endpoint
$(2, r_0)$.  At any rate, the endpoint of \eqref{strich-1} is not
essential for our regularity results.

By duality and the $TT^*$ method the estimate is equivalent to the
bilinear form estimate
$$ \int_{|s| \lesssim 1} \int_{|t| \lesssim 1} |\langle U_j(s)^* F(s), U_j(t) G(t)\rangle|\
dt ds
\lesssim 2^{-2\frac{2k-j}{q}} \|F\|_{L^\qp_t X^{\rp,2}_k}  
\|G\|_{L^\qp_t X^{\rp,2}_k}.$$
From the Hardy-Littlewood inequality
$$ \int \int |t-s|^{-\frac{2}{q}} f(t) g(s)\ dt ds \lesssim
\|f\|_\qp \|g\|_\qp,$$
valid for $q > 2$, we see that it suffices to show that
$$  |\langle U_j(s)^* F(s), U_j(t) G(t)\rangle| \lesssim
2^{-2\frac{2k-j}{q}} |t-s|^{-\frac{2}{q}} \|F(s)\|_{X^{\rp,2}_k} \|G(t)\|_{X^{\rp,2}_k}.$$
But this estimate is true for $q=\infty$, $r=2$ by the energy estimate
\eqref{energy} and Cauchy-Schwarz, while for $q=4/(n-1)$, $r=\infty$ the
result follows from the decay estimate \eqref{decay} and duality
(the fact that $q$ may be less than 1 is irrelevant).  The
general case then follows from interpolation and the assumption
\eqref{admissible}.
\end{proof}

We are almost ready to state the frequency-localized two-sided
Strichartz estimates from $L^2_t L^\rpo_x$ to $L^q_t X^{r,2}_k$.
Unfortunately the optimal exponents for these estimates depend
in a complicated way on the frequency scales $j$ and $k$.
Define the convex piecewise linear function $\alpha(j,k)$
for $j, k \geq 0$ as
$$ \alpha(j,k) = \left\{ \begin{array}{ll}
\frac{2n}{n-1}k-\frac{n+1}{n-1}j & \hbox{ for } 0 \leq j \leq k\\
2k-j & \hbox{ for } k \leq j \leq 2k\\
0 & \hbox{ for } 2k \leq j\\
\end{array}\right.
$$
Equivalently, we may
define $\alpha(j,k)$ to be the largest convex function such that
\begin{equation}\label{alpha}
\alpha(0,k) = \frac{2n}{(n-1)}k, \quad
\alpha(k,k) = k, \quad
\alpha(2k,k) = 0, \quad
\alpha(k,0) = 0
\end{equation}
for all $k \geq 0$.

\begin{proposition}\label{two-sided} If $(q,r)$ is a sharp wave-admissible 
pair, $j, k$ are non-negative integers, and $u$ is 
the solution to \eqref{standard-wave},
then
$$ 2^{\frac{\alpha(j,k)}{q}} \| S_j u\|_{L^q_t X^{r,2}_k} + \| S_j u\|_{C(L^2_x)}
\lesssim \|S_j f\|_2 + 2^{-j} \|S_j g\|_2 + 
2^{-\gamma_0 j} \|S_j F\|_{L^2_t L^{\rpo}_x}.$$
\end{proposition}

Most of these estimates are proved by the existing Strichartz estimates
and the embeddings mentioned in the previous section. 
The gain occurs
when $k \leq j \leq 2k$, so that $\alpha(j,k) = 2k-j$.  One can
show using bump function examples and a combination of parallel Knapp
examples that the estimates above are sharp, but we will not do so here.

\begin{proof} The claim involving $\| S_j u\|_{C(L^2_x)}$ follows
directly from Proposition \ref{keel}, since $(2,r_0)$ is sharp
wave-admissible and
$$ 2^{-j} 2^{\frac{n+1}{2(n-1)j}} = 2^{-\gamma_0 j}.$$
Thus it remains to treat the contribution of
$2^{\frac{\alpha(j,k)}{q}} \| S_j u\|_{L^q_t X^{r,2}_k}$.
When $\alpha(j,k) = \frac{2n}{n-1}k-\frac{n+1}{n-1}j$ this follows
from Proposition \ref{keel} and the estimate
$$ \| S_j u\|_{L^q_t X^{r,2}_k} \lesssim 2^{\frac{2nk}{(n-1)q}}
\| S_j u\|_{L^q_t L^r_x}$$
which follows from \eqref{holder} and \eqref{admissible}.  

Similarly
when $\alpha(j,k) = 0$ this follows
from Proposition \ref{keel} and the estimate
$$ \| S_j u\|_{L^q_t X^{r,2}_k} \lesssim \| S_j u\|_{L^\infty_t X^{2,2}_k}
\lesssim \|S_j u\|_{C(L^2_x)},$$
which follows from H\"older's inequality, the time localization, and
the inclusion \eqref{inclusion}.  

Thus it remains to consider the case when $\alpha(j,k) = 2k-j$,
so that $k \leq j \leq 2k$.  The contribution of $u_0$ is dealt with
in Proposition \ref{one}, so to finish the argument it suffices by \eqref{pm}
to show that
\begin{equation}\label{retarded}
 2^{\frac{2k-j}{q}} 2^{-j} \| \int_{s < t} U_j(t) U_j(s)^* F(s)\ ds \|_{L^q_t
X^{r,2}} \lesssim 2^{-\gamma_0 j} \|F\|_{L^2_t L^\rpo_x}
\end{equation}
for all Schwarz functions $F$.  As is unfortunately the case in these
types of estimates, the retarded integral \eqref{retarded} requires far more
technical manipulation than the one-sided estimates proved earlier.

When $q=\infty, r=2$ \eqref{retarded} follows from Proposition \ref{keel}, so it
suffices to verify \eqref{retarded} for the endpoint $q=2, r=r_0$.  
We will adapt the argument in \cite{tao:keel}.  By
duality \eqref{retarded} now becomes
$$
\int\int_{s < t} 
|\langle U_j(s)^* F(s), U_j(t)^* G(t)\rangle|\
dt ds
\lesssim 2^{\frac{1}{2}(j-2k)} 2^{\frac{(n+1)}{2(n-1)}j}
\|F\|_{L^2_t L^\rpo_x} \|G\|_{L^2_t X^{\rpo,2}_k}.$$
It will suffice to show that
$$
\sum_i |T_i(F,G)|
\lesssim 
 2^{\frac{1}{2}(j-2k)} 2^{\frac{(n+1)}{2(n-1)}j}
\|F\|_{L^2_t L^\rpo_x} \|G\|_{L^2_t X^{\rpo,2}_k},$$
where for $i \leq 0$, $T_i(F,G)$ denotes the bilinear form
$$ T_i(F,G) = \int\int_{t-s \sim 2^i} 
|\langle U_j(s)^* F(s), U_j(t)^* G(t)\rangle|\
dt ds.$$
By Proposition \ref{interp}
it will suffice to show that
$$ 2^{\beta(a,b) i} T_i(F,G)
\lesssim 2^{\gamma(a,b,k,j)}
\|F\|_{L^2_t L^\ap_x} \|G\|_{L^2_t X^{\bp,2}_k}$$
for all $ i \in \Z$ and
$(\frac{1}{a},\frac{1}{b})$ in a neighbourhood of $(\frac{1}{r_0},
\frac{1}{r_0})$, where
$$ \beta(a,b) = \frac{n-1}{2}(\frac{2}{r_0} - \frac{1}{a} - \frac{1}{b}).$$
$$ \gamma(a,b,k,j) = \frac{n-1}{2}(\frac{1}{2}-\frac{1}{b})(j-2k) +
\frac{n+1}{2}(\frac{1}{2}-\frac{1}{a})j.$$
By localization and time translation invariance it suffices to show that
\begin{equation}\label{beta}
|T_i(F,G)| \lesssim 2^{-\beta(a,b)i} 2^{\gamma(a,b,k,l)}
\|F\|_{L^2_t L^\ap_x} \|G\|_{L^2_t X^{\bp,2}_k}
\end{equation}
whenever $F(t), G(s)$ are supported on the time interval $|t|, |s| \lesssim 
2^{i}$.
We will prove this for the exponent pairs $(a,b) = (\infty, \infty)$, 
$(2,2)$, $(r_0, 2)$, and $(2, r_0)$, since the claim then follows by interpolation
and the fact that $2 < r_0 < \infty$  (cf. \cite{tao:keel}).

To prove the estimate when $(a,b) = (\infty, \infty)$ we use
\eqref{asymmetric} and duality to obtain
$$
|\langle U_j(s)^* F(s), U_j(t)^* G(t)\rangle|
\lesssim (2^{k-\frac{n}{n-1}j} |t-s|)^{-\frac{n-1}{2}} \|F(s)\|_1
\|G(t)\|_{X^{1,2}_k}.
$$
Integrating this over $|t-s| \sim 2^i$ we obtain
$$
|T_i(F,G)| \lesssim 2^{-\frac{n-1}{2}i} 2^{-\frac{n-1}{2}k} 2^{\frac{n}{2}j}
\|F\|_{L^1_t L^1_x} \|G\|_{L^1_t X^{1,2}_k},$$
and \eqref{beta} follows from H\"older's inequality and some algebra.

Similarly, when $(a,b) = (2,2)$, we use Cauchy-Schwarz and energy estimates
to obtain
$$ |\langle U_j(s)^* F(s), U_j(t)^* G(t)\rangle| \lesssim
\|F(s)\|_2 \|G(t)\|_{X^{2,2}_k},$$
which after integration becomes
$$ |T_i(F,G)| \lesssim \|F\|_{L^1_t L^2_x} \|G\|_{L^1_t X^{2,2}_k},$$
and \eqref{beta} again follows from H\"older's inequality.

When $(a,b) = (r_0,2)$ we write
$$ |T_i(F,G)| = |\int_t \left\langle \int_{t-s \sim 2^i} U_j(s)^* F(s)\ ds,
U_j(t)^* G(t)\right\rangle\ dt|$$
and use Cauchy-Schwarz and \eqref{energy} to obtain
$$ |T_i(F,G)| \lesssim \sup_t \|\int_{t-s \sim 2^i} U_j(s)^* F(s)\ ds\|_2
 \|G\|_{L^1_t X^{2,2}_k}.$$
However, from Proposition \ref{keel} we obtain
$$ \|\int_{t-s \sim 2^i} U_j(s)^* F(s)\ ds\|_2 \lesssim
2^{\frac{(n+1)j}{2(n-1)}} \|F\|_{L^2_t L^\rpo_x},$$
and by inserting this into the previous estimate we obtain \eqref{beta}
after using H\"older's inequality.

The case $(a,b) = (2,r_0)$ is similar. Proceeding in analogy with the
previous case we have
$$ |T_i(F,G)| \lesssim \sup_s \|\int_{t-s \sim 2^i} U_j(t)^* G(t)\ dt\|_2
 \|F\|_{L^1_t L^2}.$$
But from the adjoint of \eqref{one} we obtain
$$ \|\int_{t-s \sim 2^i} U_j(t)^* G(t)\ dt\|_2 \lesssim
2^{-\frac{2k-j}{2}} \|G\|_{L^2_t X^{r_0,2}_k},$$
and inserting this into the previous estimate we obtain \eqref{beta}
after using H\"older's inequality.
\end{proof}

\section{Proof of main theorem}

Suppose that $n > 3$, $0 < \gamma < \gamma_0$, and \eqref{concentration}
and \eqref{technical} hold.  Since $\gamma < \gamma_0$
we have from \eqref{concentration} and some algebra that
$$p < p_0
= \frac{(n+1)^2}{ (n-1)^2 + 4}.$$
We also make the technical assumption that $p > \frac{n+1}{n-1}$;
the low power case $p \leq \frac{n+1}{n-1}$ can be handled by Proposition
\ref{known}, and appears in \cite{kapitanski:wayw}.
Since $n \geq 4$, our 
assumptions on $p$ thus yield
\begin{equation}\label{p-assumptions}
 \frac{n+2}{n}, \frac{n+1}{n-1} < 
 p < \frac{n+3}{n-1}, 2.
\end{equation}

Let $f,g$ be data such that
\begin{equation}\label{eps}
\|f\|_{H^\gamma} + \|g\|_{H^{\gamma - 1}} \lesssim M,
\end{equation}
for some $M > 0$.  We will show that there exists a time $0 < T \ll 1$
that depends only on $M$, $n$, $\gamma$, $p$, and the constant in
\eqref{Fconditions}, such that 
a solution $u$ to \eqref{semilinear} exists in $C(H^\gamma)$.

We write the equation \eqref{semilinear} as an integral equation
\begin{equation}\label{iter}
u = u_0 + \Box^{-1} F(u),
\end{equation}
where the notation is as in the previous section.  

By the method of Picard iteration, to show the existence of a solution
$u$ to \eqref{semilinear} it suffices to show that the map
$u \to u_0 + \Box^{-1} F(u)$ is a contraction in some metric
space that contains $u_0$.  This space will be constructed using the
numerology used to solve \eqref{local}.  
Let $r = p\rpo$, and let $q$ be defined by \eqref{admissible}.  From
\eqref{p-assumptions} we see that $(q,r)$ is sharp wave-admissible.
We also have the inequalities
\begin{align}
2 &\leq \frac{q}{p}\label{c0}\\
\frac{\gamma_0 - \gamma}{p} - \frac{2n}{(n-1)q} + \gamma + \frac{1}{q} &\geq 0\label{c1}\\
\frac{\gamma_0 - \gamma}{p} - \frac{2n}{(n-1)q} + 2\gamma &\geq 0\label{c2};
\end{align}
indeed, \eqref{c0} simplifies to $p \leq \frac{n+3}{n-1}$, while
\eqref{c1}, \eqref{c2} are equivalent to \eqref{concentration} and 
\eqref{technical} respectively.  Since we are explicitly excluding
the endpoint \eqref{yuck}, we see that at least one of \eqref{c1},
\eqref{c2} is satisfied with strict inequality.

We now iterate in the
ball
$\{u: \|u\|_* \lesssim M\}$, where the Besov-like norm
$\| \|_*$ is given by
$$ \|u\|_* = \|u\|_{C(H^\gamma)} + (\sum_{j \geq 0} \|S_j u\|_{*,j}^2)^{1/2},$$
and the partial norms $\| \|_{*,j}$ are given by
$$ \|u\|_{*,j} = 2^{\gamma j} \sup_k 2^{\frac{\alpha(j,k)}{q}} 
\| u \|_{L^q_t X^{r,2}_k}.$$

From Proposition \ref{two-sided} we see that
$$ 
\| S_j u_0\|_{*,j}
\lesssim 2^{\gamma j} (\| S_j f\|_2 + 2^{-j} \|S_j g\|_2)$$
uniformly in $j$.  Thus from \eqref{eps} and Plancherel's theorem
we thus have that
$ \| u_0\|_* \lesssim M,$
as desired.  

It remains to show that the above map is a contraction; note that
this will give existence and uniqueness in $\| \|_*$, with the
solution depending continuously on the data in $\| \|_*$,
and hence in $C(H^\gamma)$.

It suffices to show that
\begin{equation}\label{contract}
\| \Box^{-1} (F(u) - F(v)) \|_* \ll \|u-v\|_*
\end{equation}
whenever 
\begin{equation}\label{small}
\|u\|_*, \|v\|_* \lesssim M.
\end{equation}

By using Proposition \ref{two-sided}
as before, we obtain
$$ \| S_j \Box^{-1} F\|_{*,j} \lesssim 2^{(\gamma - \gamma_0) j}
\| S_j F\|_{L^2_t L^{\rpo}_x}$$
and from Proposition \ref{keel} we obtain
$$ \| S_j \Box^{-1} F\|_{C(H^\gamma)} \lesssim \left(\sum_j (2^{(\gamma - \gamma_0) j}
\| S_j F\|_{L^2_t L^{\rpo}_x})^2\right)^{1/2},$$
for all functions $F$ and $j \geq 0$.  Also, from the energy estimate
and Sobolev embedding
we have 
$$ \| P_0 \Box^{-1} F\|_{C(H^\gamma)} \sim
\| P_0 \Box^{-1} F\|_{C(L^2)} \lesssim \| P_0 F\|_{L^1_t \dot H^{-1}}
\lesssim \| P_0 F\|_{L^1_t L^{\frac{2n}{n+2}}_x}.$$
Applying all these estimates to $F(u) - F(v)$ and using Plancherel's
theorem one obtains
\begin{align*}
\| \Box^{-1} (F(u) - &F(v)) \|_* \lesssim
\|P_0 (F(u) - F(v))\|_{L^1_t L^{\frac{2n}{n+2}}_x} +\\
&\left(\sum_{j \geq 0} (2^{(\gamma - \gamma_0) j}
\| S_j (F(u) - F(v))\|_{L^2_t L^{\rpo}_x})^2 \right)^{1/2}.
\end{align*}
Thus \eqref{contract} will follow from the non-linear estimates
\begin{equation}\label{po}
\|P_0 (F(u) - F(v))\|_{L^1_t L^{\frac{2n}{n+2}}_x} \ll \| u - v \|_{C(H^\gamma)}
\end{equation}
and

\begin{equation}\label{pre-F}
\begin{split}
 \sum_{j \geq 0} (2^{(\gamma - \gamma_0) j}
 \| S_j (F(u) - &F(v))\|_{L^2_t L^{\rpo}_x})^2
\ll\\
&\| u - v \|_{C(H^\gamma)}^2 +
 \sum_{j \geq 0} \| S_j (u - v)\|_{*,j}^2.
\end{split}
\end{equation}

We first deal with the low-frequency estimate
\eqref{po}, which is very easy. 
From 
\eqref{Fconditions}, H\"older's inequality we have
$$ \| F(u) - F(v) \|_{C(L^{2/p}_x)} \lesssim \| u - v \|_{C(L^2)}
(\| u\|_{C(L^2)}^{p-1} + \|v\|_{C(L^2)}^{p-1}).$$
u
By another application of H\"older's inequality, \eqref{small},
and the inclusion $L^2 \subset H^\gamma$ we thus have
$$ \| F(u) - F(v)\|_{L^1_t L^{2/p}_x} \lesssim T M^{p-1} \| u - v \|_{C(H^\gamma)}.$$
But since $P_0$ is given by convolution with a bump function,
\eqref{po} follows from Young's inequality (if $T$ is sufficiently
small), since one has
$1 < \frac{2}{p} < \frac{2n}{n+2}$ from \eqref{p-assumptions}.

We now turn to the high-frequency estimate
\eqref{pre-F}.  We require the following estimates.

\begin{lemma}  There exists an $\eps > 0$ such that
\begin{equation}\label{uv}
\| S_k (F(u) - F(v))\|_{L^2_t L^{\rpo}_x} 
\lesssim T^\eps
2^{\frac{2nkp}{(n-1)q}}
\| u-v\|_{L^q_t X^{r,2}_k} (\| u \|_{L^q_t X^{r,2}_k}^{p-1}
+ \| v \|_{L^q_t X^{r,2}_k}^{p-1}).
\end{equation}
for any $u$, $v$.
\end{lemma}

\begin{proof}
From \eqref{c0}, H\"older's inequality, and the definition of $r$
we have
$$
\| S_k(F(u) - F(v))\|_{L^2_t L^{\rpo}_x} \lesssim T^\eps
\| S_k(F(u) - F(v))\|_{L^{q/p}_t L^{r/p}_x}$$
for some $\eps > 0$.  From Lemma \ref{smoothing} we have
$$
\| S_k ( F(u) - F(v) )\|_{L^{q/p}_t L^{r/p}_x} \lesssim
2^{nk(1 - \frac{p}{r})} \| F(u) - F(v) \|_{L^{q/p}_t X^{r/p,1}_k}.$$
But from \eqref{Fconditions} and 
H\"older's inequality we have
$$ \| F(u) - F(v) \|_{L^{q/p}_t X^{r/p,1}_k}
\lesssim \| u-v \|_{L^{q}_t X^{r,p}_k} 
( \| u \|_{L^q_t X^{r,p}_k}^{p-1} + \| v \|_{L^q_t X^{r,p}_k}^{p-1}).$$
By \eqref{holder} and \eqref{p-assumptions} the right-hand side is dominated by
$$  2^{-nk(1 - \frac{p}{2})}
\| u-v\|_{L^q_t X^{r,2}_k} 
(\| u \|_{L^q_t X^{r,2}_k}^{p-1} + \| v \|_{L^q_t X^{r,2}_k}^{p-1}).$$
Combining all these estimates and using \eqref{admissible} the
lemma follows.
\end{proof}

\begin{lemma} There exist $i \in \{1,2\}$ and $\eps > 0$ such that
\begin{equation}\label{gen-decay}
  \| f \|_{L^q_t X^{r,2}_k} \lesssim
  2^{\frac{\gamma_0 - \gamma}{p} k} 2^{-\frac{2nk}{(n-1)q}}
  ( 2^{-\eps k} \| f\|_{C(H^\gamma)} + \sum_j 2^{-\eps|j-ik|} 
  \| S_j f\|_{*,j} ).
\end{equation}
for all $f$.
\end{lemma}

\begin{proof}
From \eqref{triangle} it suffices to show that
\begin{equation}\label{po-est}
 \| P_0 f\|_{L^q_t X^{r,2}_k} \lesssim
2^{\frac{\gamma_0 - \gamma}{p} k} 2^{-\frac{2nk}{(n-1)q}}
2^{-\eps k} \| f\|_{C(L^2)}
\end{equation}
and
\begin{equation}\label{sj-est}
\|S_j f\|_{L^q_t X^{r,2}_k}
\lesssim
  2^{-\frac{\gamma - \gamma_0}{p} k} 2^{-\frac{2nk}{(n-1)q}}
  2^{-\eps|j-ik|}
    \| S_j f\|_{*,j}.
\end{equation}
We first consider the low-frequency estimate \eqref{po-est}.  From
Proposition \ref{young} and H\"older's inequality we have
$$ \| P_0 f\|_{L^q_t X^{r,2}_k} \lesssim 2^{-nk(\frac{1}{2} - \frac{1}{r})}
\| f\|_{C(L^2)},$$
and so \eqref{po-est} reduces to
$$ -n(\frac{1}{2} - \frac{1}{r}) \leq \frac{\gamma_0 - \gamma}{p} 
-\frac{2n}{(n-1)q} - \eps,$$
which reduces using \eqref{admissible} to the hypothesis
$\gamma < \gamma_0$.

We now turn to \eqref{sj-est}.  From
the definition of $\| \|_{*,j}$ it suffices to 
show that
$$ 1 \leq 2^{-\frac{\gamma - \gamma_0}{p} k}
2^{-\frac{2nk}{(n-1)q}}
    2^{-\eps|j-ik|} 2^{\gamma j} 2^{\frac{\alpha(j,k)}{q}}.$$
uniformly in $j$ and $k$, for some $\eps > 0$.  This reduces to showing that
$$ (\frac{\gamma_0 - \gamma}{p} - \frac{2n}{(n-1)q})k
+ \gamma j + \frac{\alpha(j,k)}{q} \geq \eps |j-ik|.$$
By the convexity of $\alpha$ it suffices to verify this inequality
for the four ranges in \eqref{alpha}.  Dividing by $k$,
it thus suffices to verify that
\begin{align*}
(\frac{\gamma_0 - \gamma}{p} - \frac{2n}{(n-1)q})
+ \frac{2n}{(n-1)q} &\geq \eps |0-i|\\
(\frac{\gamma_0 - \gamma}{p} - \frac{2n}{(n-1)q})
+ \gamma + \frac{1}{q} &\geq \eps |1 - i|\\
 (\frac{\gamma_0 - \gamma}{p} - \frac{2n}{(n-1)q}) 
+ 2 \gamma &\geq \eps |2 - i|\\
 \gamma
&\geq \eps |1-0i|.
\end{align*}
The first and fourth inequality follow from the hypothesis 
$0 < \gamma < \gamma_0$.
From \eqref{c1} and \eqref{c2} we see that the second and third inequalities
are satisfied with $\eps = 0$, and since at least one of these inequalities
is assumed to hold with strict inequality one can make $\eps > 0$
by choosing $i$ appropriately.
\end{proof}

Applying \eqref{gen-decay} to $u$,$v$ we obtain
\begin{equation}\label{u-est}
  \| u \|_{L^q_t X^{r,2}_k},
  \| v \|_{L^q_t X^{r,2}_k}
  \leq T^\eps M 2^{-\frac{\gamma - \gamma_0}{p} k} 2^{-\frac{2nk}{(n-1)q}},
\end{equation}
since we have from \eqref{small} that
$$
2^{-\eps k} \| u\|_{C(H^\gamma)} + \sum_j 2^{-\eps|j-ik|} \| S_j u \|_{*,j}
\lesssim \| u\|_* \lesssim M,$$
and similarly for $v$.
If we also apply \eqref{gen-decay} to $u-v$, and insert the resulting
inequality and \eqref{u-est} into \eqref{uv}, one obtains
\begin{align*}
\| S_k (F(u) - &F(v))\|_{L^2_t L^{\rpo}_x}
\leq T^{\eps p} M^p
2^{-(\gamma - \gamma_0) k} \\
&(2^{-\eps k} \| u-v\|_{C(H^\gamma)} + 
\sum_j 2^{-\eps|j-ik|} \| S_j (u - v) \|_{*,j}).
\end{align*}
Thus, the left-hand side of \eqref{pre-F} is majorized by
$$
T^{2\eps p} M^{2p} \left(\| u - v \|_{C(H^\gamma)}^2 + \sum_{k \geq 0}
(\sum_{j \geq 0} 2^{-\eps|j-ik|} \| S_j u - S_j v\|_{*,j})^2
\right).$$
By the Cauchy-Schwarz inequality this is majorized by
$$ T^{2\eps p} M^{2p} \left(\| u - v \|_{C(H^\gamma)}^2 + 
\sum_{j \geq 0} \sum_{k \geq 0} 
2^{-\frac{1}{2}\eps|j-ik|} \| S_j u - S_j v\|_{*,j}^2
\right),$$
and \eqref{pre-F} follows by evaluating the $k$-summation if $T$ is 
sufficiently small.  This concludes the proof.


\begin{thebibliography}{10}

\bibitem{beals:xsb}
M. Beals, \emph{Self-Spreading and strength of Singularities for solutions to
semilinear wave equations}, Annals of Math \textbf{118} (1983), 187--214.

\bibitem{beals:bezard}
M. Beals, M. Bezard, \emph{Low regularity local solutions for
field equations}, Comm. Partial Differential Equations \textbf{21} (1996): 79--124.
\bibitem{bergh:interp}
J. Bergh, J. L\"ofstr\"om, \emph{Interpolation Spaces: An Introduction},
Springer-Verlag, 1976.

\bibitem{borg:schrodinger}
J. Bourgain, \emph{A remark on Schrodinger operators},
Israel J. Math. 77 (1992), 1--16.

\bibitem{borg:xsb}
J. Bourgain, \emph{Fourier restriction phenomena for certain lattice
subsets and applciations to nonlinear evolution equations}, Part I,
Geometric and Funct. Anal. \textbf{3} (1993), 107-156.

\bibitem{damiano}
D. Foschi, \emph{Lecture Notes for S. Klainerman's Graduate Course In 
Nonlinear Wave Equations:  Fall 1996, Princeton University}, 
Private Communication.

\bibitem{ginebre:summarywave}
J. Ginebre, G. Velo, \emph{Generalized Strichartz Inequalities for the
Wave Equation}, Jour. Func. Anal., \textbf{133} (1995), 50--68.

\bibitem{kapitanski:variablestrichartz}
L. Kapitanski, \emph{Some Generalizations of the Strichartz-Brenner
Inequality}, Leningrad Math. J., \textbf{1} (1990), 693--676.

\bibitem{kapitanski:wayw}
L. Kapitanski, \emph{Weak and Yet Weaker Solutions of Semilinear Wave
Equations}, Comm. Part. Diff. Eq., \textbf{19} (1994), 1629--1676.

\bibitem{kpv:xsb}
C.~E. Kenig, G. Ponce, L. Vega, \emph{The Cauchy problem for the Korteweg-de
Vries equation in Sobolev spaces of negative indices}, Duke Math J.
\textbf{71} (1993), 1--21.

\bibitem{klainerman:xsb}
S. Klainerman, M. Machedon, \emph{Smoothing estimates for null forms 
and application}, Duke math J.
\textbf{81} (1995), 99--103.

\bibitem{tao:keel}
M. Keel, T. Tao, \emph{Endpoint Strichartz Estimates}, to appear, Amer. Math. J.

\bibitem{lindblad:sharpduke}
H. Lindblad, \emph{A Sharp Counterexample to Local Existence of Low Regularity
Solutions to Nonlinear Wave Equations},
Duke Math J., \textbf{72}, (1993), 503--539.

\bibitem{lindbladsogge:semilinear}
H. Lindblad, C.~D. Sogge, \emph{On Existence and Scattering with Minimal
Regularity for Semilinear Wave Equations}, Jour. Func. Anal., \textbf{130} (1995),
357--426.

\bibitem{mss:localsmoothing}
G. Mockenhaupt, A. Seeger, C.~D. Sogge,  \emph{Local Smoothing of
Fourier Integrals and Carleson-Sj\"olin Estimates},
J. Amer. Math. Soc., \textbf{6} (1993), 65--130.
 
\bibitem{vargas:restrict}
A. Moyua, A. Vargas, L. Vega, \emph{Schr\"odinger Maximal Function
and Restriction Properties of the Fourier transform}, International Math.
Research Notices \textbf{16} (1996).

\bibitem{sadosky:interp}
C. Sadosky, \emph{Interpolation of Operators and Singular Integrals},
Marcel Dekker Inc., 1976.

\bibitem{siegel:strichartz}
I.~E. Segal, \emph{Space-time Decay for Solutions of Wave Equations},
Adv. Math., \textbf{22} (1976), 304--311.
 
\bibitem{sogge:smoothing}
C.~D. Sogge, \emph{Propogation of singularities and maximal functions in the
plane}, Invent. Math. \textbf{104} (1991), 349--376.

\bibitem{sogge:fio}
C.~D. Sogge, \emph{Fourier integrals in classical analysis},
Cambridge University Press, 1993.

\bibitem{sogge:nlwbook}
C.~D. Sogge, \emph{Lectures on Nonlinear Wave Equations},
International Press, 1995.
 
\bibitem{stein:small}
E.~M. Stein, \emph{Singular Integrals and Differentiability
Properties of Functions}, Princeton University Press, 1970.

\bibitem{strichartz:restrictionquadratic}
R.~S. Strichartz, \emph{Restriction of Fourier Transform to Quadratic Surfaces
and Decay of Solutions of Wave Equations}, Duke Math. J., \textbf{44} (1977), 705--774.

\bibitem{tomas:restrict}
P. Tomas, \emph{A restriction theorem for the Fourier transform}, Bull. Amer. 
Math. Soc. \textbf{81} (1975), 477--478.

\bibitem{triebel:interp}
H. Triebel, \emph{Interpolation Theory, Function Spaces, Differential
Operators}, North-Holland, 1978.
 
\end{thebibliography}
\end{document}